\documentclass{amsart}
\usepackage[mathscr]{eucal}
\usepackage{verbatim}

\title[ ]{Positive combinations of projections in von Neumann algebras and purely infinite simple
 C*-algebras  }

\author{Victor Kaftal}

\address{Department of Mathematics\\
University of Cincinnati\\
P. O. Box 210025\\
Cincinnati, OH\\
45221-0025\\
USA}

\email{victor.kaftal@math.uc.edu}

\author{P.  W. Ng}

\address{Department of Mathematics\\
University of Louisiana\\
217 Maxim D. Doucet Hall\\
P.O. Box 41010\\
Lafayette, Louisiana\\
70504-1010\\
USA}

\email{png@louisiana.edu}

\author{Shuang Zhang}

\address{Department of Mathematics\\
University of Cincinnati\\
P.O. Box 210025\\
Cincinnati, OH\\
45221-0025\\
USA}

\email{zhangs@math.uc.edu}

\newtheorem{thm}{Theorem}[section]
\newtheorem{df}{Definition}[section]
\newtheorem{lem}[thm]{Lemma}

\newtheorem{prop}[thm]{Proposition}
\newtheorem{rem}[thm]{Remark}
\newtheorem{qu}[thm]{Question}

\newcommand{\B}{\mathscr{B}}

\newcommand{\A}{\mathscr{A}}

\newcommand{\Cu}{\mathscr{O}}

\newcommand{\h}{\mathscr H}

\newcommand{\K}{\mathscr K}

\newcommand{\Mul}{\mathscr M}

\DeclareMathOperator{\tr}{Tr}
\newcommand{\BH}{\mathbb{B}(\mathscr H)}

\def\sideremark#1{\ifvmode\leavevmode\fi\vadjust{\vbox to0pt{\vss
\hbox to 0pt{\hskip\hsize\hskip1em
\vbox{\hsize2cm\tiny\raggedright\pretolerance10000
\noindent#1\hfill}\hss}\vbox to8pt{\vfil}\vss}}}

\def \bib(#1;#2;#3;#4;#5;#6)  {{#1},{\it #2} {#3},
{\bf#4} (#5) {#6}\par\smallskip}

\begin{document}

\begin{abstract}
We give an overview of the question: which positive elements in an operator
 algebra can be written as a linear combination of projections with positive coefficients.
 A special case of independent interest is the question of which positive elements can
 be written as a sum of finitely many projections. We  focus on von Neumann algebras,
  on  purely infinite simple C*-algebras, and on their  associated multiplier algebras.
\end{abstract}

\maketitle
 \section {\bf Introduction and preliminaries on projections in  C*-algebras}\label{S:Intro}

In this article we focus on a number of questions involving projections in C*-algebras.

\begin{qu}\label{Q:pos} Which positive elements in a C*-algebra $\A$
are linear combination of projections with positive coefficients? (\textit{Positive combinations of projections} for short.) And for which algebras are all positive elements positive combinations of projections?
\end{qu}

A case of particular interest is when all the coefficients are equal to one.

\begin{qu}\label{Q: finite}
Which positive elements in a C*-algebra $\A$ are finite sums of
projections?
\end{qu}

A related and somewhat simpler question for those C*-algebras where it is natural to consider also
 infinite sums (e.g., von Neumann algebras, with convergence in the strong topology or
 multiplier algebras of a C*-algebra with convergence in the strict topology) is the following:

\begin{qu}\label{Q: strong}
Which positive operators in a C*-algebra $\A$ are finite or infinite sums of projections?
\end{qu}

We will study these questions for von Neumann algebras, for
$\sigma$-unital purely infinite simple C*-algebras and for their
associated multiplier algebras. In the future,
we plan to extend our
analysis to a wider class of algebras.

Our paper is organized as follows. In Section 2, we review some of the basic
notions related to projections in C*-algebras.

In Section \ref {S:lin comb} we survey known results related to
linear combination of projections. We introduce formally an algebra
constant $v_0$ that bounds the minimal sum of the moduli of the
coefficients (see Definition \ref  {V0}). To the best of our
knowledge, this constant was first introduced for $\BH$ by Fong and
Murphy in \cite {FongMurphy}. Although often not explicitly referred
to, this constant is implicit in much of the work on linear
combinations of projections.  By combining the work of Fack \cite {Fack}
and Marcoux  \cite{MarcouxSpan} we show the existence of that
constant for a large class of C*-algebras (see Proposition \ref
{prop:universalconstant}). The reason that we are interested in the
existence of $v_0$ is that together with an additional mild
condition, it guarantees that positive \textit{invertible} elements
are positive combinations of projections (see Proposition
\ref{Positive invertible}). 

In Section \ref {S:positive} we review the results about positive
combinations of projections in $\BH$  and then present a simple but
useful lemma (Lemma \ref {2x2lemma}) that is the key to decomposing
positive elements into positive combinations of projections. We
apply it to three kinds of C*-algebras.

For von Neumann algebras, it permits us to obtain a sufficient
condition in terms of spectral resolution. This condition is also
necessary if the range projection of the element is properly
infinite and this condition is always satisfied for finite factors and
$\sigma-$finite type III factors (Theorem \ref {pos comb W*}).  For
factors, this condition is equivalent to the non-vanishing of the
essential norm with respect to the relative compact ideal (i.e., the
ideal generated by the finite projections of the factor).

For global von Neumann algebras, a similar condition (see Theorem \ref {T3.4})
 holds by replacing essential norm with \textit{central essential norm} and
 invoking the theory of central ideals of Halpern \cite {Hh72}, \cite {Hh77}
  and Stratila and Zsido \cite {SsZl73} (see Section \ref{S:positive} for these notions).

One might argue that $\sigma$-unital purely infinite simple
C*-algebras share some properties with $\sigma$-finite type III
factors.  The analogy holds in this case: all positive elements in
either type of algebras are positive combinations of projections
(Theorems \ref {pos comb W*}(iii),
 \ref  {pos comb in C*}).
Moreover, the same holds for multiplier algebras of $\sigma$-unital purely
 infinite simple C*-algebras (Theorem \ref {multipl positive}).

In Section \ref {S:infinite}, we discuss decompositions into infinite
sums of projections. Their characterization is complete for $\BH$
and  for type III $\sigma-$finite factors, while for type II factors
we characterize diagonal elements that are  infinite sums of
projections.  (Theorem \ref {strong sum W}). The characterization is
also complete  for multiplier algebras of $\sigma$-unital purely
infinite simple C*-algebras (Theorem \ref {mult infin}), or more
precisely for those elements in $\Mul(\A)\setminus \A$.

Finally, in Section \ref {S:finite}, we discuss decompositions into finite sums of
projections. A complete characterization  is still open in $\BH$. We review results
 in this area, including the sufficient condition $\|x\|_e>1$ in terms of the
 essential norm, initially considered by Choi and Wu \cite{WpCm09}, (see Theorem \ref {Choi-Wu}).

We extend this condition to von Neumann factors (Theorem \ref {T:finsumsfactors})
and to global von Neumann algebras in terms of the central essential spectrum
(Theorem \ref {T4.3}).

Ideas from frame theory lead to a necessary condition in terms of operator ideals
 of $\BH$. This shows that the full characterization of finite sums of projections
 may require some delicate analysis.

The sufficient condition $\|x\|_e>1$  for $x$ to be a finite sum of
projections holds also in multiplier algebras of $\sigma$-unital
purely infinite simple C*-algebras (Theorem \ref {multipl fin
sums}). It would thus be natural to think  that positive elements of
$\sigma$-unital purely infinite simple C*-algebras with operator
norm greater than $1$ may  also be finite sums of projections.
Indeed, this turns out to be the case  for the Cuntz algebras $\Cu
_n$ with $2\le n< \infty$  (Theorem \ref {Cuntz}),  but at the moment we do not
know if the property
 holds or not for $\Cu _\infty$.

\section{Some facts about projections and C*-algebras}\label{S:projections}

During the last five decades, much research activity around
operator algebras  has centered at  the structure of projections and
their ``levels of abundance". In some C*-algebras, the amount of
projections can be abundant enough to give rise to an approximate
spectral decomposition of all normal elements, while some other
C*-algebras do not have any nontrivial projections besides the zero
projection or the identity (even in certain simple C*-algebras).

\medskip
 Let us first recall several notions concerning the ``amount" of projections in a C*-algebra.

\begin{df} A C*-algebra $\ A$  can have the following properties:

\item[($\exists $P)]   There exists a nontrivial projection in $\A
$ (\cite{Blackadar}).

\item[(AP)] $\A $ is generated as a C*-algebra by its projections (\cite{Blackadar}).

 \item[(SP)]  There exists a nontrivial projection in each nonzero  hereditary
  C*-subalgebra of $\A $  (\cite{Blackadar}).

 \item[(LP)]  $\A $ is the closure of the span of its projections  (\cite{Blackadar},\cite{ZhangSimple}).
 Equivalently, the
 self-adjoint part  $\A _{sa}$ of $\A $ is the
closure of the collection of linear combinations with real
coefficients of projections in  $\A $.

In fact, if (LP) holds  and $a=a^*\in \A $, then $a= \lim \sum_1^n
\lambda_jp_j$, and hence, $a= a^*=\lim   (\sum_1^n \lambda_jp_j
)^*$. It follows that $a = \lim \sum_1^n \frac{\lambda_j+ \bar
\lambda _i}2 p_j$. The converse is also trivial.

\item [(HLP)]  Every hereditary C*-subalgebra of $\A $ has  the LP property
(\cite{KNZ2010}).

\item[(CP)] The positive part  $\A _+ $ of $\A$ is the closure of all linear
 combinations with positive coefficients of projections  of
$\A$  (\cite{Blackadar}).

\item[(PLP)]  Every hereditary C*-subalgebra of $\A $ has the property
CP (\cite{KNZ2010}).

 \item[(HP)] Every hereditary C*-subalgebra of $\A $ has an approximate  identity
 consisting of projections (\cite{Blackadar}).

 \item[(FS)] Every self-adjoint element of $\A $ can be approximated in
 norm by self-adjoint elements with finite spectra. It is equivalent
 to replace ``self-adjoint" by ``positive" (\cite{Blackadar}).

 \item[(RR=0)] The collection of self-adjoint invertible elements is norm-dense in
 the  self-adjoint part of $\A $ (\cite{BP}).

\item[(AUD)] $\A $ is called almost upward directed, if for any
positive number $\epsilon $ and any two projections $p, q$ in $\A $
there exists another projection $r \in \A $ in the hereditary C*-subalgebra generated
by $p, q$  such that $\|(1-r)p\|<
\epsilon $ and $\| (1-r) q\| <\epsilon $ (\cite{Zhang3}).
\end{df}

We also recall the following facts.

\begin{rem}
\item[(1)] The three properties HP, FS, and RR=0 are equivalent and the algebras possessing them are called of  \textit{real rank zero}. These are the algebras showing the most abundance of projections
 (\cite{Blackadar} and \cite{BP}). And the property  AUD for a multiplier
algebra $M(\A)$ is equivalent to $RR(M(\A )) =0$ if $\A $ is
 $\sigma$-unital and of real rank zero (\cite{Zhang3}).
\item [(2)] All hereditary C*-subalgebras of a C*-algebra of real rank zero, including all von Neumann algebras,  have the property LP   (\cite{BP}, \cite{Blackadar}).
\item [(3)] All hereditary C*-subalgebras of the multiplier algebra of a  stable $\sigma$-unital C*-algebra
 of real rank zero have the property LP (\cite{ZhangSimple}).
\end{rem}

In this paper, we focus on three kinds of operator algebras: von Neumann algebras, $\sigma$-unital purely
infinite simple C*-algebras, and their multiplier algebras.

Among von Neumann algebras, of main interest are of course the
$\sigma$-finite factor algebras;
 roughly speaking, we can think of
algebras operating on a separable Hilbert space. Then the (infinite)
type I factors  are of course identified with the algebra $\BH$ of
bounded operators on a Hilbert space $\h$, the type III factors have
no proper ideals nor finite nonzero projections, and the type
II$_\infty$ factors can be viewed as a continuous analogs of $\BH$
endowed with the ideal of ``relative compact" operators.

A simple C*-algebra is called purely infinite if every hereditary
C*-subalgebra of the form $(x\A x^*)^-$ for some $x\in \A$ contains
an infinite projection $p$ in the sense of Murray-von Neumann,
namely, $p\sim p'<p$.  It is easy to see that  all hereditary
C*-subalgebras of a separable C*-algebra have that form.  Familiar
examples of such C*-algebras include all type III factors,  and the
Cuntz algebras $\Cu _n$ $(n\le 2\le \infty )$ constructed around
1980, which were the first separable examples of purely infinite
simple C*-algebras.

Around 1980 both Cuntz and Blackadar were wondering about ``What the
abundant level of projections  is in $\Cu _n$s". This question is
clear from the following item (3) in terms of real rank zero.

\begin{rem}
\item[(1)] All von Neumann algebras have real rank zero (\cite{BP}).

\item[(2)] All purely infinite simple C*-algebras have real rank zero
(\cite{Zhang1} and \cite{Zhang4}).

\item [(3)] A simple C*-algebra $\A $ is purely infinite  if and only if $\A $
has real rank zero and every nonzero projection is infinite (\cite{Zhang2} and
\cite{Zhang7}). Notice the analogy between these algebras and type
III factors.

\item[(4)] A $\sigma$-unital purely infinite simple $C^*$-algebra is either
unital  or it is stable, i.e. $\A \cong \A \otimes \K$
(\cite{Zhang4}).

\end{rem}

Concerning topologies of the multiplier algebras, we recall that
besides the one induced by the norm, $\Mul(\A)$ is endowed with the
\textit{strict} topology defined by: $ x_{\lambda } \longrightarrow
x$ in the strict topology if and only if
$$\text{max}\{\|(x_\lambda -x)a\|, \|(x_\lambda ^*-x^*)a\| \}
\longrightarrow 0\ \ \ \text{for all}\ \ a\in\A.$$ If $\A =\K $,
then $\Mul(\A ) =\BH$. In this extreme case, the strict convergence
is exactly  *-strong convergence (also called the double strong
convergence). Thus it  is natural  to consider strictly converging
series of projections in $\Mul(\A )$.

The following properties illustrate the richness of the structure of
$\Mul(\A )$.

 \begin{rem}
 \item[(a)] Assume
that $\A $ is simple. Then the generalized Calkin algebra $\Mul (\A
\otimes \K)/{\A \otimes \K }$ is simple if and only if either $\A $
is elementary (i.e., finite matrices or $\K$) or  $\A $ is purely
infinite (\cite{ZhangSimple}).

\item[(b)]  If $\A $ is $\sigma$-unital  purely infinite and simple, then every
projection in $\Mul(\A ) \setminus  \A $ is equivalent to the
identity. (A close analogue of $\B(\h )$) (\cite{ZhangSimple}).

\item[(c)] If $\A $ is $\sigma$-unital purely infinite simple, then the
generalized Calkin algebra $\Mul(\A )/\A $ is purely infinite and
simple (\cite{ZhangSimple}).

\item[(d)]  Assume that $\A $ is $\sigma$-unital, purely infinite and  simple.
Then $\Mul(\A )$ has real rank zero if and only if $K_1(\A )=\{0\}$
(\cite{Zhang4} and \cite{ZhangSimple}).
\end{rem}

\section{\bf Linear combinations of projections   }\label{S:lin comb}

The structure of projections and their level of
abundance have always been of central interest in the study of operator algebras. On the one hand,
projections in a von Neumann algebra are so abundant to permit the
spectral decomposition of all normal elements, and on the other
hand, some operator algebras, even some simple unital C*-algebras,
do not have any nontrivial projection at all.

In this article, we focus on algebras where every element is a linear
combination of projections.  It was first proven in 1967 by Fillmore
that this holds for the algebra $\mathbb{B}(\h)$ of all bounded
operators on a separable Hilbert space $\h$.

\begin{thm}\label{T:Fillmore1967}(\cite{Fillmore2})
  If $\h$ is an
infinite dimensional separable Hilbert space, then every operator in
$\mathbb{B}(\h)$ is a linear combination of at most $257$
projections.
\end{thm}
This estimate was subsequently improved by Pearcy and
Topping \cite{PearTop} in the same year.

\begin{thm}\label{T:PearcyTopping1967}(\cite{PearTop})
If $\h$ is an infinite dimensional separable Hilbert space, then
every operator in $\mathbb{B}(\h)$ is a linear combination of at
most   $16$ projections.
\end{thm}

In 1984 Matsumoto \cite{Matsumoto}  showed the following best
estimate in $\mathbb{B}(\h)$:

\begin{thm}\label{T:Matsumoto1984}(\cite{Matsumoto})
If $\h$ is an infinite dimensional separable Hilbert space, then
every selfadjoint operator in $\mathbb{B}(\h)$ is a linear
combination of $5$ projections. Consequently, every operator in
$\mathbb{B}(\h)$ is a linear combination of $10$ or fewer
projections.
\end{thm}

Pearcy and Topping generalized the $\mathbb{B}(\h)$ result  to infinite
von Neumann factors:

\begin{thm}\label{T:PearcyTopping1969}(\cite{PearTop})
  All elements in an infinite von Neumann factor are  linear combinations  of projections.
 \end{thm}

The more difficult case of type $II_1$ factors  was settled by
Fack and de La Harpe.

\begin{thm}\label{T:FackHarpe1980}(\cite{FackHarpe})
  All elements in a $II_1$ factor are  linear combinations of projections.
  \end{thm}

 Thus all elements in a von Neumann factor are  linear combinations of projections.
 Beyond factors, this property can fail. In 1992 Goldstein and Paszkiewicz  proved   the following:

\begin{thm}\label{T:GP1992}(\cite[Theorems 3.2]{GP})
A von Neumann algebra   has the property that all its elements are  linear combinations
of projections if and only if the algebra has no finite type I direct  summand
with infinite-dimensional center. \end{thm}

In their work, Goldstein and Paszkiewicz improved the previous estimates on the number
 of needed projections. However, in the present article, we will not focus on those estimates,
 but rather on the size of the coefficients needed in the linear combinations, more precisely,
 on the existence of an algebra constant   $v_0$ with the property:

\begin{df}\label{V0}
A C*-algebra $\A$ has a constant $v_0$ if every $x\in \A$  can be
decomposed into
 a linear combination $x= \sum_{j=1}^{n}t_jp_j$ for certain number of projections $p_j\in
 \A$ (here $n$ depends on $x$) such that
$
  \sum_{j=1}^{n}|t_j|\le v_o\|x\|.
$
\end{df}
 Such a universal constant $v_0$ exists for every von Neumann algebra with no  type I direct  summand
  with infinite-dimensional
 center; the fact is implicitly given  in the statements of the Goldstein and Paszkiewicz's
 results that are restated as follows.

\begin{thm}\label{T:W*-v0}(\cite[Theorems 1--3]{GP})
 If a von Neumann algebra $M$ has no finite type I direct  summand with infinite-dimensional
  center, then it has a constant $v_o$ as in Definition \ref {V0}. \\  If $M$
  is properly infinite then $v_o\le8;$\\
if $M$ is of type II$_1$, then $v_o\le14;$\\
if $M$ is the direct sum of $m$ matrix algebras, then  $v_o\le m+4$.
\end{thm}

\bigskip

The situation is of course more complicated for  C*-algebras. There
is a large class of  C*-algebras including all unital simple
AF-algebras with finitely many extreme tracial states and all unital
purely infinite simple C*-algebras for which every element is equal
to a linear combination of projections. We also know of the
following results:

\begin{thm}\label{T:MarcouxMurphy1998} The following C*-algebras $\A$
have the property that every element is a  linear combination of projections in the algebra:
\begin{enumerate}
\item  $\A$ is a unital simple real rank zero
C*-algebra with unique tracial state $\tau$ and satisfying strict comparison
for projections: if $p, q \in \A$ are projections and $\tau(p) < \tau(q)$ then  $p \prec q$
  (\cite{MarcouxCommutator}).
\item $\A$ is a unital simple C*-algebra with a nontrivial projection
and does not admit a tracial state (\cite{MarcouxCommutator},
\cite{MarcouxMurphy}).
\item $\A$ is a unital simple AT-algebra with real rank zero and finitely
many extreme tracial states (\cite{MarcouxSpan}, \cite{MarcouxCommutator}).
\item $\A$ is a unital simple AH-algebra with real rank zero, bounded
dimension growth and finitely many extreme tracial states such that
the building blocks are finite direct sums of full matrix algebras over
unital commutative C*-algebras (with uniformly bounded dimension)
(\cite{MarcouxSpan}, \cite{MarcouxCommutator}).
\end{enumerate}
\end{thm}

Using the work of  \cite {Fack} and \cite{MarcouxSpan}   we obtained   the following

\begin{prop}\label{prop:universalconstant}(\cite [Proposition 2.6] {KNZ2010})
Let $\A$ be a unital $C^*$-algebra for which there exist two
projections $p, q\in \A$ such that $p \sim q \sim 1 $ and $pq=0$.
Then there exists a universal constant $v_0$ as in Definition \ref{V0}
\end{prop}

This class of C*-algebras includes among others all unital simple infinite
C*-algebras and all multiplier algebras of $\sigma$-unital
stable C*-algebras.

\begin{qu}
If a C*-algebra $\A$ has the property that every element is a linear
combination of projections in the algebra, must there exist a
constant $v_o$ as in Definition \ref{V0}?
\end{qu}

\section{\bf Positive combinations of projections }\label{S:positive}
If in a C*-algebra all elements are  linear combinations of
projections,
 it is obvious that self-adjoint elements can be expressed as linear combinations
  of projections with real coefficients.  What about expressing positive elements
  as linear combinations of projections with positive coefficients? (\textit{positive combinations
  of projections} for short).

\

Fillmore remarked in \cite {Fillmore2} that positive compact
operators of infinite rank cannot be decomposed as  positive
combinations of projections. Later in 1985, Fong and Murphy proved in
\cite{FongMurphy} that these operators are the only exceptions.

\begin{thm}\label{T:FongMurphy1985}(\cite{FongMurphy})
A positive operator $x\in \mathbb{B}( \h)$ can be written as a
positive combination of projections  if and only if $x$ is not a
compact operator of infinite rank.
\end{thm}

The key for obtain this result was the fact that positive \textit{invertible} operators
 are positive combinations of projections. This fact was obtained in \cite {Fillmore}
 by Fillmore  from his characterization of sums of two projections (see Theorem \ref{2 proj})
  and again by  Fong and Murphy in \cite{FongMurphy} by the existence of a constant $v_o$ as
   in Definition \ref {V0}. The  Fong and Murphy proof  can be modified  to prove the same result for
   certain C*-algebras.

\begin{prop}\label{Positive invertible}\cite[Proposition 2.7] {KNZ2010}
Let $\A$ be a unital C*-algebra with  the following two properties:
\begin{enumerate}
\item  $\A$ has a constant $v_0$ as in Definition \ref {V0};
\item  $\A$ satisfies condition CP, namely,  positive combinations of projections are norm dense in $\A_+$.
\end{enumerate}
Then every positive invertible element in $\A$ is a positive combination of projections.
\end{prop}
The class of C*-algebras satisfying the above conditions include all
von Neumann algebras with no finite type I direct summand with
infinite dimensional center, all unital simple purely infinite
C*-algebras, and all multiplier algebras of stable $\sigma$-unital
 C*-algebras of real rank zero.

The following lemma is the key to extend the positive combination property beyond invertible elements.
\begin{lem} \label{2x2lemma}(\cite[Lemma 2.9]{KNZ2010})

Let $\A$ be C*-algebra with two mutually orthogonal projections
$q\prec p$ with the property that
 for every projection $r\in \A$ such that
$r\le p$, every positive invertible element of $r\A r$ is a
positive combination of projections in  $  r\A r$. Let $x:=t p+ b$
be a positive element  in $\A$, where $t
>\|b\|$ and $b=qb=bq$.   Then $x$
is a positive combination of projections in $\A$.
\end{lem}

The application of this lemma to von Neumann algebras yields the following characterization of
 positive combinations of projections in terms of their spectral resolution.

 \begin{thm}\label{pos comb W*}(\cite [Theorem 2.12]  {finitesums})
 Let $M$ be a von Neumann algebra,  $x \in M^+$, $r_x$ be the range projection of $x$, and $\chi(x)$
 be the spectral measure of $x$.
 \item[(i)] Assume that $r_xMr_x$ has no finite type I direct summands with infinite dimensional
 center and that    there is a $ \delta > 0$  such that  $\chi_{(0, \delta)}(x)\prec \chi_{[\delta,\infty)}(x)$.
 Then $x$ is a positive combination of projections in $M$.
 \item[(ii)] If $r_x$ is properly infinite, then $x$ is a positive combination of projections in
 $M$ if and only if  there is a $ \delta > 0$  such that  $\chi_{(0, \delta)}(x)\prec \chi_{[\delta,\infty)}(x)$.
  \item[(iii)]  If  $r_xMr_x$  is a finite sum of finite factors or of  $\sigma$-finite type III factors,
  then $x$ is always a positive combination of projections in $M$.
 \end{thm}

For finite algebras the spectral condition $\chi_{(0, \delta)}(x)\prec \chi_{[\delta,\infty)}(x)$ for some $\delta >0$  in Theorem \ref {pos
comb W*}(i)  is only sufficient.

\begin{qu}
Are there  natural necessary and sufficient conditions that characterize those elements
that are positive combinations of projections in a finite von Neumann algebra with infinite
dimensional center?
\end{qu}

Notice that if $M=\BH$ for   $\h$ separable and if $x$ is an
infinite rank positive
 operator, then the condition $\chi_{(0, \delta)}(x)\prec \chi_{[\delta,\infty)}(x)$ for some
  $\delta>0$ holds if and only if  $x$ is not compact, and the latter condition
  is indeed  the necessary and sufficient condition in Theorem \ref {T:FongMurphy1985}. The same
   equivalence holds in type II$_\infty$ factors where  in lieu of the ideal of compact
    operators we take the ideal $J(M)$ of \textit{relative compact} operators,
   which has been studied first by Breuer \cite {Bm68}, \cite {Bm69} and Sonis \cite{Sm71}, i.e.,
   the norm closed two sided-ideal generated by the finite projections of $M$. Denote by
    $\|\cdot\|_e$ the essential norm corresponding to $J(M)$, i.e.,   $\|x\|_e=\|\pi(x)\|$, where
     $\pi: M\to M/ J(M)$ is the quotient map.

Thus the factor case is well understood and simple: either  $x$
belongs to a finite factor, in which case $x$ is  always a positive
combination of projections or $r_x$ is infinite, in which case $x$
is a positive combination of projections if and only if $x\not \in
J(M)$, or, equivalently, $\|x\|_e>0$.
 If $M$ is of type III, then $ J(M)=\{0\}$ and again $x$ is  always a positive combination of projections.

In order to move beyond factors, or finite sums of factors, we need
to reformulate   the condition $\chi_{(0, \delta)}(x)\prec \chi_{[\delta,\infty)}(x)$ for some $\delta>0$  in terms of an
appropriate ideal. The natural tool to do so is the theory of
central ideals of Halpern \cite {Hh72}, \cite {Hh77} and Stratila
and Zsido \cite {SsZl73}, of which we will briefly review here only
a couple of key notions.

A two-sided norm-closed  ideal $\mathscr J$ of $M$ is called
\textit{central} if for every collection of mutually orthogonal
projections $e_\gamma\in M\cap M'$ and every norm-bounded family
$x_\gamma \in \mathscr J$ it follows that $\sum_\gamma x_\gamma
e_\gamma \in \mathscr J$. Given a central ideal $\mathscr J$ and an
$x\in M$ there is a \textit{central essential spectrum} of $x$
relative to $\mathscr J$, which is a bounded and strongly closed
subset of the center. If $x\ge 0$, that subset has a maximal element
$C_{u, \mathscr J}(x)$ which we call the \textit{central essential
norm} of $x$.

For instance, if $M=\bigoplus_{k=1}^\infty \mathbb B(\h_k)$ with
$\h_k\equiv \h$,  $\mathscr J=   \bigoplus_{k=1}^\infty K(\h_k)$ and
$x=\bigoplus_{k=1}^\infty x_k\in M^+$, then $\mathscr J$ is a
central ideal of $M$ and $C_{u, \mathscr J}(x)=
\bigoplus_{k=1}^\infty \|x_k\|_eI_k$ where $\|x_k\|_e$ is the
essential norm of $x_k$.

Given  a properly infinite projection $p\in M$, there is an
associated central ideal  $\mathscr J(p)$ of ``dimension less than
$p$", namely  the ideal whose set of projections is
$$\{q\in M \mid  q\le c(p) \text{ and }pg\prec  qg \text{ for some projection }
 g\in M\cap M' \Rightarrow pg=0.\}$$
where $c(p)$ denotes the  central support of $p$.

While in an infinite factor the necessary and sufficient
 condition for a  positive element $x$ to be a positive combination of
 projections is that its essential norm is nonzero, for global algebras
 we must require that the central essential norm be \textit{bounded away from zero}.
 More precisely:

\begin{thm}\label{T3.4}(\cite [Theorem 3.4]  {finitesums})
Let $M$ be a von Neumann algebra, let $x \in M^+$ and assume that
$r_x$ is properly infinite. Set $\mathscr J= \mathscr J( r_x)$. Then
$T$ is a positive combination of projections in $M$ if and only if
 $C_{u,\mathscr J} (x)$ is locally invertible, i.e., there is a $\nu > 0$ for which
 $C_{u,\mathscr J} (x)\ge \nu c(r_x)$.
\end{thm}

It is often easier in properly infinite von Neumann algebras  to
  decompose  elements into linear combinations
of projections,
  positive combinations of projections, or sums of projections, and and to
  investigate commutators  than in the finite
von Neuamnn algeras.

Thus in moving now our focus to C*-algebras, we naturally start with
$\sigma$-unital purely infinite simple C*-algebras  and  their
associated multiplier algebras. Both these kinds of algebras satisfy
the conditions for  Proposition \ref  {Positive invertible} and
Lemma \ref {2x2lemma} and Lemma \ref {2x2lemma}  plays  again a key
role.

Perhaps not surprisingly the  situation for $\sigma$-unital purely infinite simple
C*-algebras is  similar to the one for type III factors:

\begin{thm} \label{pos comb in C*}\cite [Theorem 2.11]{KNZ2010} Let $\A $ be a $\sigma$-unital
 purely infinite simple C*-algebra.
Then every nonzero  positive  element of $\A $  is a positive
combination of projections.  \end{thm}

The crux of the proof is the case when the hereditary algebra
generated by $x$ is unital and the number $\|t\|$  is not an
isolated point in the spectrum of $x$ and then the proof consists of  
a reduction to the situation in Lemma \ref {2x2lemma}. In the
general case, the hereditary algebra generated by $x$ can be
embedded in a unital hereditary subalgebra of $\A$ and then the
result also follows.

A similar proof combined with results for $\BH$ from \cite{WpCm09} and \cite{finitesums}, shows that if $x\in\Mul(\A )^+\setminus
\A$ then $x$ is a scalar multiple of a finite sum of projections. Here, the
analogy is to the case when $x\in \BH^+\setminus K(\h)$. But then,
combining this result with the one of Theorem \ref {pos comb in C*},
we get:\begin{thm}\label{multipl positive} (\cite[Theorem 3.1]
{KNZ2010}) Let $\A $ be a $\sigma$-unital  purely infinite simple
C*-algebra. Then every nonzero  positive  element of $\Mul(\A )$  is
a positive combination of projections.
\end{thm}

\section{\bf Infinite sums of projections }\label {S:infinite}

A motivation for studying infinite sums of projections in $\BH$ with
convergence in the strong topology \textit{strong sums} for short,
comes from frame theory. Indeed Question \ref{Q: strong} for $\BH$
is reformulation of the frame theory question of which positive operators are the Bessel operators of
a uniform-norm Bessel sequence. With this motivation, Dykema,
Freeman, Kornelson, Larson, Ordower and Weber in \cite{Larsonetal}
investigated Question \ref{Q: strong} and found that a sufficient
condition for a positive operator $x$ to be a sum of projections is
that  the essential norm $\|x\|_e>1.$ A more general question -- 'which
positive operators are the Bessel operators of a Bessel sequence
with pre-assigned norms' -- was studied by Kornelson and Larson in 2004
(\cite {Larson1}), and by Antezana, Massey, Ruiz and Stojanoff in
2007 (\cite{AMRS}), leading to an alternative approach to the
Schur-Horn Theorem.

 A full answer to Question \ref{Q: strong} was recently obtained by the authors
 of the present article for the case when $M$ is a $\sigma$-finite type I or type III
  von Neumann factor,  and for diagonalizable operators in the case when $M$ is a
  $\sigma$-finite type II factor. (An operator $x\in M^+$ is  diagonalizable if
   $x =\sum \gamma_j e_j$ for some  $\gamma_j>0$  and  mutually orthogonal projections $e_j\in M$.)
   The characterization of strong sums of projections is given in terms of the notions of
   \textit{excess part} and \textit{defect part} of an operator defined as follows:
\begin{align*} &\chi & &\text{the characteristic function }\\
&r_x=\chi_{(0, ||x||]}(x) & &\text{the range projection of $x$}\\
&x_+:=(x-I)\chi_{(1, ||x||]}(x) & & \text{the excess part of}\ \ x \\
&x_-:= (I-x)\chi_{(0,1)}(x) & & \text{the defect part of } \ \ x.
\end{align*}  Thus we have the decomposition $$x= x_+ - x_-+r_x.$$

Notice that $\|x\|_e>1$ if and only if $x_+$ is not compact.

\begin{thm} \label{strong sum W}(\cite[Theorem 1.1]{B(H)case})
\item [(i)] Let  $M=\BH$ and $x\in M^+$. Then $x$ is a strong sum of projections
if and only if either $\tr(x_+)=\infty$ or $\tr(x_-)\le
\tr(x_+)<\infty$  and $\tr(x_+) - \tr(x_-)\in \mathbb N\cup \{0\}$.
\item[(ii)] Let  $M$ be a $\sigma$-finite  type II factor, $\tau$ be a faithful,
semifinite, normal trace on $M^+$, and $x\in M^+$.  If $x$ is a
strong sum of projections, then $\tau(x_+)\ge \tau(x_-)$. If either
$\tau(x_+)=\infty$ or  $x$ is diagonalizable and $\tau(x_+)\ge
\tau(x_-)$, then T is  a strong sum of projections.
 \item[(iii)] Let $M$  be a $\sigma$-finite  type III factor and $x\in M^+$. Then
  $x$ is a strong sum of projections if and only if either $||x||>1$ or $x$ is a projection.
 \end{thm}
The proof is based on a $2\times2$ matrix construction that
decomposes  linear combinations $(1-\lambda)p_1+(1+\mu)p_2$ of two
equivalent projections $p_1$ and $p_2$, with $\lambda,\mu >0$,
$\lambda\le 1$,  into the combination $q_1+ (1+\mu-\lambda)q_2$ of
two equivalent projections $q_1$ and $q_2$.  A reduction to strong
sums of such linear combinations is always possible in the type I or
type III case, or in the type II case when $\tau(x_+)=\infty$. To
achieve this reduction in the type II case when $\tau(x_+)<\infty$,
we ask for $x$ to be diagonalizable.  However, it is easy to find
examples of strong sums of projections in type II$_1$ factors, even
sums of two projections, that are not diagonalizable.  This leads us
to ask:
   \begin{qu}
 Can the hypothesis that $x$ is diagonalizable be removed from (ii)?
 \end{qu}

In many respects, the strict topology in the multiplier algebra $\Mul(\A)$ of a
 C*-algebra $\A$  takes the role of the strong topology in von Neumann factors.
 Multiplier algebras have also a natural notion of essential norm, namely  for every
  $x\in \Mul(\A)$,  we denote by $\|x\|_{e}$ the norm of $\pi(x)$ where
$\pi:\Mul(\A)\to \Mul(\A)/ \A$ is the canonical quotient map. When
$\A$ is sufficiently ``nice", which in this context is
$\sigma$-unital, nonunital, purely infinite simple, we have the
following characterization of those elements that are infinite sums
of projections.

 \begin{thm} \label{mult infin} (\cite[Theorem 1.1 ]{multipliercase}) Let $\A $ be a
$\sigma$-unital, nonunital, purely infinite simple  $C^*$-algebra
and $x$ be a positive element of $\Mul(\A )$. Then $x$ is a strictly
converging sum of projections belonging to $\A$ if and only if one
of the following mutually exclusive conditions hold:
\item [(i)] $\| x\|_{e} > 1$.
\item [(ii)] $\| x \|_{e} = 1$ and $\|x \| > 1$.
\item [(iii)] $x\in \Mul (\A )\setminus \A $ is a projection.
\item[(iv)]  $x$ is the sum of finitely many projections
belonging to $\A$. \end{thm}

A key step in the proof of this theorem, whose crux is case (ii), and a result
of independent interest in itself, is the following:

\begin{lem} (\cite[Lemma 2.5] {multipliercase}) For every
$\sigma$-unital, nonunital, purely infinite simple  $C^*$-algebra
$\A $ the collection of  finite sums of projections is norm dense in
$\{x\in \A^+\mid \|x\|>1\}$.
\end{lem}

Notice that case (iv) in Theorem \ref {mult infin} leads naturally to the harder Question
 \ref {Q: finite} (see the next section).

\section{\bf Finite sums of projections }\label{S:finite}
The question of which positive operators are finite sums of projections, i.e.,
Question \ref{Q: finite}, is harder, and  indeed it is still open in $\BH$.

 The first results in this direction were obtained in  1969 by Fillmore:

\begin{thm}\label {fin rank} \cite {Fillmore}
A finite rank positive operator $x\in \mathbb{B}(\h)$ is a sum of
projections if and only if $\tr (x)\in \mathbb N$ and $\tr (x) \ge
$rank($x$).
\end{thm}

\begin{thm}\label {2 proj}\cite {Fillmore}
A positive operator $x\in \mathbb{B}(\h)$ is a sum of two
projections if and only if $x= 2p\oplus b$  where $p$ is a possibly
zero projection and $b$ is a positive operator  unitarily equivalent
to $2I-b$.
\end{thm}

In 1988, Choi and Wu  obtained the following result that was announced
in a survey article \cite [Theorem 4.12] {Wpy94} by Wu   and was presented together
with other related results  in their recent paper \cite[Theorem 2.2] {WpCm09}.

\begin{thm}\label {Choi-Wu}\cite {WpCm09}
A sufficient condition for a positive operator $x\in   \BH$ to be a
finite sum of projections is that the essential norm $\|x\|_e>1$.
\end{thm}

A special case of this result, i.e., for $x=\alpha I$  was obtained
through completely different methods  by Kruglyak,  Rabanovich, and
Samo\u{i}lenko (see \cite {KRS02}, \cite {KRS03},  and several other
articles referenced there) who characterized the values of $\alpha
>1$ for which $\alpha I$ is the sum of $n$ or more projections.

The authors together with Halpern recently re-obtained independently the Choi
and Wu  sufficient condition  in the more general setting of von Neumann algebras.
 For $\sigma-$finite factors the formulation is very similar to the $\BH$ one.

\begin{thm}\label{T:finsumsfactors}\cite [Corollary 4.4 ] {finitesums}  Let $M$
be an infinite $\sigma-$finite factor and $x\in M^+$. A sufficient
condition for $x$ to be a finite sum of projections in $M$ is that
\item [(i)] $\|x\|_e>1$ when $M$ is of type I$_\infty$ (usual essential norm of $\BH$);
\item [(ii)] $\|x\|_e>1$ when $M$ is of type II$_\infty$  (essential norm relative to the ideal
  $ J(M)$ generated by the finite projections of $M$);
\item [(iii)]  $\|x\|>1$ when $M$ is of type III (operator norm).
\end{thm}

The  proof is based on decomposing the operator into a sum of finite rank
summands where the range projections of non-consecutive summands are orthogonal,
 and where each summand is a finite sum of projections, with the number of projections
  needed for each summand being uniformly bounded.

Naturally, for global properly infinite  von Neumann algebras the essential norm has
 to be replaced by the \textit{central essential norm} $C_{u,\mathscr J}(x)$ with respect
 to the central ideal $\mathscr J(r_x)$ with ``dimension less than" $r_x$ (see Section 4.)

\begin{thm}\label {T4.3}\cite [Theorem 4.3] {finitesums}  Let $M$ be a von Neumann algebra. Assume that
  $x\in M^+$ has a properly infinite range projection $r_x$
and $C_{u,\mathscr J}(x) \ge \nu c(r_x)$ for some $\nu>1$. Then $x$
is a finite sum of projections in $M$. \end{thm}

For the case of  type II$_1$ factors we  have a sufficient condition in
terms of the excess and defect part of the operator, albeit again only for
the case of diagonalizable operators, (cf. Theorem \ref {strong sum W}).

\begin{thm}\label{T6.5}\cite[Theorem 6.5]   {finitesums}
Let $M$ be a  type II$_1$ factor and let $x\in M$ be a positive
diagonalizable operator. Then a sufficient condition for  $x$ to be
a finite sum of projections in $M$ is  that  $\tau(x_+) >
\tau((x_-).$
\end{thm}

It is easy to see  that the condition that the essential norm
$\|x\|_e>1$ is not necessary for $x$ to be a finite sum of
projections. In fact, a necessary condition, which for simplicity's
sake we formulate just for $\BH$, is

\begin{thm}\label{nec cond}(\cite[Corollary 5.7] {finitesums}
Let $x\in \BH^+$ be a finite sum of projections and assume that
$x_+\in K(H)$. Then $x_-\in K(H)$. If $x_-=0$ then $x_+$ has finite
rank. If $x_-\ne 0$ then $x_+$ and $x_-$ generate the same two-sided
ideal of $\BH$.
\end{thm}

This necessary condition permits to easily construct examples of
positive operators (with essential norm $=1$),  that are strong
sums of projections but cannot be finite sums of projections.

In passing now to C*-algebras, we notice that the same proof that
permits to decompose a positive element $x\in  \Mul(\A )\setminus \A
$ into a positive combination of projections (see Theorem \ref
{multipl positive}) actually guarantees that if $\|x\|_e>1$, then
one of the coefficients $t_j>1$ and thus reduce the problem to the
case that $x=tp\oplus sq$ where $t>1$ and $s\ge 0$.  By considering
$p$ and $q$ as infinite projections in a copy of $\BH$ embedded
unitally in a corner of $\Mul(\A)$, this reduces the problem to
$\BH$:

 \begin{thm}\label{multipl fin sums}(\cite[Theorem 3.3]  {
 KNZ2010})   \label{thm:pismultiplierN} Let $\A $ be a $\sigma$-unital
 but nonunital purely infinite simple C*-algebra.
Then every  positive element  $x\in \Mul(\A )\setminus \A $ with
$\|x\|_{e}>1$ is a finite sum of projections.
 Consequently, every nonzero positive element of $\Mul(\A )\setminus \A $ can
be written as a positive scalar multiple of a sum of projections.
  \end{thm}

The problem when  $\|x\|_{e}=0$, i.e.,  $x\in \A$ seems harder.  So
far we have obtained  this result which did surprise us:

\begin{thm}\label{Cuntz} If $2\le n<\infty $, then every
positive element $x$ in the Cuntz algebra  $\Cu _n$ with $\|x\|>1$
is a  sum of finitely many projections. \end{thm}

\section*{Acknowledgements}
S. Zhang would like to acknowledge that his research was partially
supported by Taft International Travel Grant. P. W. Ng would like to 
acknowledge that his research was done while visiting the 
other two authors at the University of Cincinnati, in trips supported
by the UC Math Department.

\end{document}